# THE MEAN, VARIANCE AND LIMITING DISTRIBUTION OF TWO STATISTICS SENSITIVE TO PHYLOGENETIC TREE BALANCE

By Michael G. B. Blum, Olivier François and Svante Janson

*University Joseph Fourier, CNRS, Grenoble, Institut National Polytechnique*
*de Grenoble and Uppsala University*

For two decades, the Colless index has been the most frequently used statistic for assessing the balance of phylogenetic trees. In this article, this statistic is studied under the Yule and uniform model of phylogenetic trees. The main tool of analysis is a coupling argument with another well-known index called the Sackin statistic. Asymptotics for the mean, variance and covariance of these two statistics are obtained, as well as their limiting joint distribution for large phylogenies. Under the Yule model, the limiting distribution arises as a solution of a functional fixed point equation. Under the uniform model, the limiting distribution is the Airy distribution. The cornerstone of this study is the fact that the probabilistic models for phylogenetic trees are strongly related to the random permutation and the Catalan models for binary search trees.

**1. Introduction.** *Phylogenetic trees* (PT) represent the shared history of extant species. The idea of using trees to model evolution dates back to Darwin ([10], see his diagram on page 117). In a (rooted) PT, there is a common ancestral species called the *root* and each branching represents the time at which a divergence has occured. A PT is usually reconstructed using data from $n$ different species (or *taxa*) which are located at the *leaves*. The tree has $n-1$ internal nodes that correspond to the ancestors of the sample. There are two distinct features of rooted PTs. First is the branching structure or topology of the tree. Second is the branch lengths which correspond to periods of time separating major evolutionary events. The shape of such trees carries useful information about the history of diversification rates among species by reflecting the footprint left by evolutionary processes.

Biologists have extensively investigated the ways in which the shapes of PTs can be measured [20]. Mooers and Heard [25] wrote an exhaustive review concerning tree balance in systematic biology and Aldous [5] gave an introduction in a more mathematical setting. How these measures are related to macroevolution processes has been studied by Rogers [31] and Agapow and Purvis [1], relying upon intensive computer simulations. Thus far, several statistics have been introduced to measure the shapes of PTs (see Agapow and Purvis [1], for eight of them). Among these statistics, the most widespread are Sackin's and the Colless indices. *Sackin's index* [34, 36] counts the number of ancestors crossed in the path from each leaf to the root. *Colless' index* [9] looks at the internal nodes, partitioning the leaves that descend from them into groups of sizes $L$ and $R$, and computes the sum of absolute values $|L - R|$ for all ancestors.

The probability distributions of the Sackin and the Colless statistics have been investigated for various models of biologically plausible random trees. Two random models of PT are often considered in the literature. The most famous is the *Yule model* [38]. The Yule model is a branching process with constant speciation rate where the number of extant species is specified. The assumption of a constant speciation rate may be weakened by assuming that the diversification rate could vary in time but is the same for all species at any time. This assumption does not modify the distribution of PT shape. An alternative model considered by biologists is called the *uniform model*. It assumes that all PTs are equally likely. This model is biologically motivated, as it arises from a large family of Galton–Watson processes conditioned by the total size of the trees (see [2]). In addition, McKenzie and Steel [24] have shown that when speciation events are constrained to occur before a time $\tau$ after their previous speciation event, the resulting process converges to the uniform model as $\tau$ tends to zero. In both models, Rogers [31] studied the joint distribution of the Sackin and the Colless statistics using numerical computations. He concluded that these statistics were strongly correlated in large PTs. The limiting distribution of the Colless statistic was also conjectured to be non-Gaussian [30].

This article describes the mean, the (co)variance and the limiting joint distribution of Sackin's and Colless' indices for large PTs under the Yule and uniform models. Because this study is mainly concerned with the topology of PTs, branch lengths can be ignored. A PT is then a cycle-free connected graph with vertices of degree one (the leaves), two (the root) or three (all ancestors except the root). Leaves are usually labeled, whereas ancestors are not. This simplified model of phylogeny without branch lengths is sometimes called a *cladogram* (see [4]). Our proofs use the connection to recent results in theoretical computer science, as well as the correspondence between PTs and random binary search trees (BST). This approach extends results by Blum and François [6] who showed that Sackin's index has the same limit distribution as the number of comparisons used by the quicksort algorithm [16]. More specifically, we deal with the Yule and uniform models separately. For the Yule model, our analysis relies on the recursive structure of the tree and makes use of the fixed-point method (see, e.g., [17]). This method was introduced in the probabilistic analysis of algorithms by Rösler [32]. In the uniform model, the results are based on the connection between uniform trees and Bernoulli excursions [37]. A large family of statistics similar to Sackin's and Colless' indices have been studied by Fill and Kapur [12] under the Catalan model for BSTs.

In Section 2, we shall present our main results. Section 3 explains how probabilistic models for PTs are related to probabilistic models for BSTs. Section 4 is dedicated to the Yule model, while Section 5 deals with the uniform model.

**2. The Sackin and the Colless statistics.** Consider a PT with $n$ leaves. The Sackin statistic adds the number of internal nodes between each leaf and the root of the tree to produce the index

$$S_n = \sum_{i=1}^{n} d_i,$$

where the sum runs over the $n$ leaves of the tree and $d_i$ is the number of ancestors crossed in the path from $i$ to the root (including the root). The Colless statistic looks at the internal nodes, partitioning the leaves that descend from them into groups of sizes $L_j$ and $R_j$ and computing

$$C_n = \sum_{j=1}^{n-1} |L_j - R_j|,$$

where the sum runs over the internal nodes and $L_j$ (resp. $R_j$) corresponds to the number of leaves in the left (resp. right) subtree under node $j$.

Denote by $\mathcal{M}_2$ the space of all bivariate, centered probability measures with finite second moments, and by $\mathcal{L}(X)$ an element of $\mathcal{M}_2$. We have the following result:

THEOREM 1. *Assume the Yule model of PT. Consider the map $\mathcal{T}: \mathcal{M}_2 \to \mathcal{M}_2$ such that for all $\nu \in \mathcal{M}_2$, we have*

$$\mathcal{T}(\nu) = \mathcal{L}\left(\begin{bmatrix} U & 0 \\ 0 & U \end{bmatrix} \begin{pmatrix} S \\ C \end{pmatrix} + \begin{bmatrix} 1-U & 0 \\ 0 & 1-U \end{bmatrix} \begin{pmatrix} S' \\ C' \end{pmatrix} + \begin{pmatrix} b_S \\ b_C \end{pmatrix}\right),$$

*with*

$$\begin{pmatrix} b_S \\ b_C \end{pmatrix} = \begin{pmatrix} 2U\log U + 2(1-U)\log(1-U) + 1 \\ U\log U + (1-U)\log(1-U) + 1 - 2\min(U, 1-U) \end{pmatrix},$$

where $(S,C)$, $(S',C')$ and $U$ are independent random variables such that $\mathcal{L}(S,C) = \mathcal{L}(S',C') = \nu$ and $U$ is uniform over the interval $(0,1)$. Then we have

$$\left(\frac{S_n - \mathbb{E}[S_n]}{n}, \frac{C_n - \mathbb{E}[C_n]}{n}\right) \xrightarrow{d} (S,C), \qquad n \to \infty,$$

where the convergence holds in distribution and the limiting probability distribution is the unique fixed point of the map $\mathcal{T}$.

REMARK 1. The convergence in Theorem 1 will actually be proven for a stronger topology than convergence in distribution. As can be seen from Section 4, it indeed holds for the Wasserstein–Mallows $d_2$-metric [29] which guarantees the existence and convergence of the second moments.

REMARK 2. This result extends the fact that the normalized Sackin index

$$\bar{S}_n = \frac{S_n - \mathbb{E}[S_n]}{n} \tag{1}$$

converges in distribution to the same limit as the number of comparisons in the quicksort algorithm. According to Rösler [32], the limit $S$ satisfies a (functional) fixed-point equation of the type

$$S \stackrel{d}{=} US + (1-U)S' + 2U \log U + 2(1-U)\log(1-U) + 1, \tag{2}$$

where $S, S'$ and $U$ are independent random variables, $S$ and $S'$ are identically distributed, $U$ is uniformly distributed over the inverval $(0,1)$ and the identity holds for distributions. Regarding Colless' index, the functional fixed-point equation becomes

$$\begin{aligned}C \stackrel{d}{=} &UC + (1-U)C' + U \log U \\ &+ (1-U)\log(1-U) + 1 - 2\min(U, 1-U).\end{aligned} \tag{3}$$

A well-known result in systematic biology is that the expectation of $S_n$ is of order $2n \log n$. More precisely, Kirkpatrick and Slatkin [20] showed that

$$\mathbb{E}[S_n] = 2n \sum_{j=2}^{n} \frac{1}{j}$$

and

$$\mathbb{E}[S_n] = 2n \ln n + (2\gamma - 2)n + o(n),$$

where $\gamma$ is Euler's constant. Using the connection to the quicksort algorithm, the variance of the limiting distribution can be obtained according to Knuth [21] as

$$\text{Var}[S_n] \sim \left(7 - 2\frac{\pi^2}{3}\right)n^2, \qquad n \to \infty.$$

These results can be extended to the case of Colless index as follows, taking into account that $S_n$ and $C_n$ are strongly correlated for large PTs:

THEOREM 2. *Assume the Yule model of PT. Then we have*

(4) $$\mathbb{E}[C_n] = n\log n + (\gamma - 1 - \log 2)n + o(n),$$

(5) $$\text{Var}[C_n] \sim \left(3 - \frac{\pi^2}{6} - \log 2\right)n^2,$$

(6) $$\text{Cor}[S_n, C_n] \sim \frac{27 - 2\pi^2 - 6\log 2}{\sqrt{2(18 - \pi^2 - 6\log 2)(21 - 2\pi^2)}} \approx 0.98,$$

*as $n$ goes to infinity.*

Regarding the uniform model of PT, mathematical results have received less attention than for the Yule model. After an appropriate rescaling, we prove the convergence of both $S_n$ and $C_n$ to the same marginal probability distribution and identify this distribution as $\sqrt{8}$ times the integral of the standard Brownian excursion $e(t)$,

$$\omega = \int_0^1 e(t)\,dt.$$

The distribution of random variable $\mathcal{A} = \sqrt{8}\omega$ is known as the *Airy distribution*. A formula for the moments of $\mathcal{A}$ has been given by Flajolet and Louchard [13]. In particular, we have

$$\mathbb{E}[\mathcal{A}] = \sqrt{\pi}$$

and

$$\text{Var}[\mathcal{A}] = \frac{10 - 3\pi}{3}.$$

THEOREM 3. *Assume the uniform model of PT. Then we have*

(7) $$\frac{S_n - C_n}{n^{3/2}} \xrightarrow{p} 0$$

*and*

$$\frac{S_n}{n^{3/2}} \xrightarrow{d} \mathcal{A},$$

*as $n$ goes to infinity.*

REMARK 3. Regarding $S_n$, the connection to the internal path length of a BST enables us to immediately state that

$$\frac{S_n}{n^{3/2}} \xrightarrow{d} \mathcal{A}.$$

This result was actually established by Takacs [37] using the method of moments. In addition, we find that

$$\mathbb{E}[S_n] \sim \sqrt{\pi} n^{3/2}$$

and

$$\mathrm{Var}[S_n] \sim \left(\frac{10}{3} - \pi\right) n^3.$$

The moments of $C_n$ follow from the next theorem.

THEOREM 4. *Assume the uniform model of PT. Then we have*

$$\mathbb{E}[C_n] \sim \sqrt{\pi} n^{3/2}$$

*and*

$$\mathrm{Var}[C_n] \sim \frac{10 - 3\pi}{3} n^3$$

*as $n$ goes to infinity. In addition, the variables $S_n$ and $C_n$ are asymptotically correlated, that is,*

$$\mathrm{Cor}[S_n, C_n] \sim 1$$

*and we have, for any $k, \ell \geq 0$,*

$$\mathbb{E}[C_n^k S_n^\ell] \sim n^{3(k+\ell)/2} \mathbb{E}[\mathcal{A}^{k+\ell}]$$

*as $n$ goes to infinity.*

REMARK 4. While $C_n$ and $S_n$ are by far the most popular statistics used in studies of phylogenetic imbalance, other measures have also been considered (see [1]). Some of these can also be studied in the Yule model using the contraction method, mainly because they are defined as sums of elementary functions of subtrees over all nodes. For example, the result for the Fusco and Conk statistic, modified by Purvis, Katzouralus and Agapov [28], is left to the reader. In the same spirit, we believe that the $B_1$ index of Shao and Sokal [36] could be studied without difficulties. Studying the remaining statistics ($B_2$ and $\sigma_N^2$) would nevertheless require considerably more effort.

REMARK 5. In a recent large-scale study of the phylogenetic database, Blum and François [7] considered the shape statistic

$$F_n = \sum_{j=1}^{n-1} \log(N_j - 1),$$

where the sum runs over all internal nodes and $N_j$ represents the number of extant descendants of internal node $j$. A similar statistic had been previously proposed by Chan and Moore [8], but the logarithm was omitted. Once the normalizing constant has been removed, $F_n$ corresponds to the logarithm of the probability of a tree in the Yule model. In particular, the statistical test based on $F_n$ is the most powerful test for rejecting the Yule model against the uniform model and conversely (Neyman–Pearson theorem). Fill [11] showed that $F_n$ has a Gaussian distribution (for large trees) and gave asymptotic expansions for the means and variances under both the Yule and the uniform models (see also [12]).

**3. Phylogenetic and binary search trees.** Trees are often encountered in theoretical computer science as data structures associated with *divide and conquer algorithms*. In this section, we explain how binary search trees can be mapped onto phylogenetic trees univoquely and how probabilistic models for BSTs are transported on probabilistic models for PTs.

*Mapping binary search trees.* A binary tree can be defined recursively. It is either empty or it is a node (the root) with left and right subtrees. A binary search tree is a binary tree where labels are associated with the vertices. These labels are constrained: the label of a vertex is greater than or equal to all labels contained in the left subtree and less than or equal than all labels contained in the right subtree (Section 5.5, [35]). The transformation that maps BSTs into PTs can be found in [4]. Given a BST with $n-1$ vertices, the structure is modified as follows. Vertices in the BST become ancestors in the PT. To accomplish this, two leaves are connected to each vertex of degree one and one leaf is connected to each vertex of degree two. The root has a special status. If the degree is 0, 1 or 2, then 2, 1 or 0 leaves are added. The labels of leaves are chosen arbitrarily from the $n!$ possible orders.

Two obtained PTs are equivalent if their left and right subtrees can be interchanged recursively (see Figure 1). The set of PTs is the set of equivalence classes for this equivalence relation. Figure 2 gives a graphical representation of these transformations.

A PT may therefore arise from the construction of $2^{n-1}$ equivalent modified BSTs. Because there are $\mathcal{C}_{n-1}$ BSTs with $n-1$ vertices, we obtain $\mathcal{C}_{n-1} n!/2^{n-1}$ possible PTs. This number coincides with the total number of PTs, which equals $(2n-3)!!$, where

$$(2n-3)!! = (2n-3)(2n-5)\cdots 3 \cdot 1.$$

*Probabilistic models.* The mapping described in the above paragraph also transfers probabilistic models for BSTs to probabilistic models for PTs. For instance, there will be an equivalence between the random permutation model for BSTs [22] and the Yule model for PTs. Probabilistic models for BSTs (with $n-1$ vertices) can be described as a general class of models called *branching Markov processes*. A definition of branching Markov processes can be found in [4]. We recall this definition here. Let $\hat{q}_{n-1}$ be a symmetric probability distribution on $\{0, \ldots, n-2\}$:

$$\hat{q}_{n-1}(i) = \hat{q}_{n-1}(n-2-i), \qquad i = 0, \ldots, n-2.$$

In the branching Markov process, the size of the left subtree is chosen according to the probability distribution $\hat{q}_{n-1}$. This procedure is repeated recursively in subtrees, assuming local independence. The probability distribution $\hat{q}_{n-1}$ is called the *splitting distribution*. In the same way, probability distributions on PTs with $n$ leaves can be associated with splitting probability distributions $q_n$ on $\{1,\ldots,n-1\}$. At each step, the labels of the left subtree, of size $i$, are sampled uniformly from the $\binom{n}{i}$ possible labels. At the end of the construction, left–right distinctions are simply suppressed in the building of the PT.

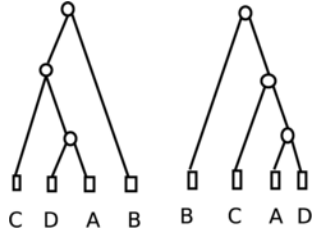

FIG. 1. *Two graphical representations of the same PT. They are seen to be identical by interchanging left and right subtrees.*

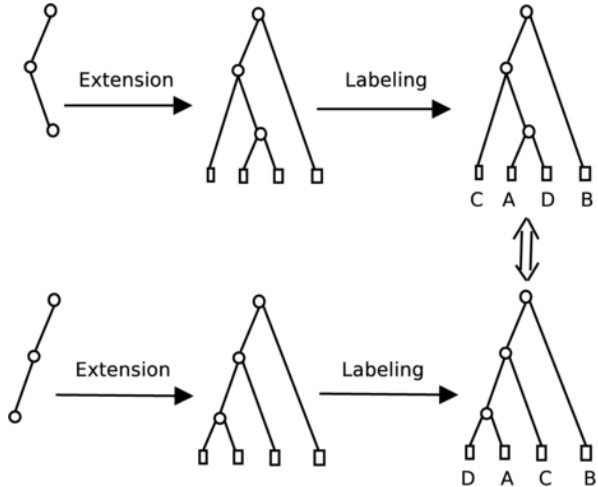

FIG. 2. *The transformation of two binary search trees to the same PT. The extension consists of connecting two leaves to vertices with outdegree* 0 *and one leaf to vertices with outdegree* 1. *The two resulting trees represent the same PT.*

LEMMA 1. *Assume that $\hat{T}_{n-1}$ is a BST, sampled according to a branching Markov process with splitting probability $\hat{q}_{n-1}$. Denote by $\Phi_n$ the transformation which consists of extending a BST with $n-1$ vertices to a PT with $n$ leaves. Then $T_n = \Phi_n(\hat{T}_{n-1})$ is a PT, sampled according to a branching Markov process with splitting probability $q_n$ such that*

$$q_n(i) = \hat{q}_{n-1}(i-1), \qquad i=1,\ldots,n-1.$$

PROOF. This is a consequence of the basic properties of $\Phi_n$. If a BST has $i$ vertices in its left subtree, the resulting PT has $i+1$ leaves in one of the two subtrees of the root. The symmetry property of $\hat{q}_n$ ensures that all members of the same equivalence class have the same probability of occurrence. □

Lemma 1 has the interesting consequence that well-studied models of BSTs can be transposed into models on PTs. The Yule and uniform models for PTs then get associated with special cases of branching Markov processes.

On one hand, the random permutation model for BSTs with $n-1$ vertices is a branching Markov process with splitting probability

$$\hat{q}_{n-1}(i) = \frac{1}{n-1}, \qquad i=0,\ldots,n-2.$$

This model is mapped by $\Phi_n$ into the Yule model for PTs with $n$ leaves with splitting probability

$$q_n(i) = \frac{1}{n-1}, \qquad i=1,\ldots,n-1.$$

The splitting distribution for Yule trees was found by Harding [15]. Note that the same splitting property also holds for $n$-coalescent tree topologies [19].

On the other hand, the Catalan model for binary BSTs with $n-1$ vertices assumes that all $\mathcal{C}_{n-1}$ binary trees have the same probability of occurrence. The number of trees with a left subtree of size $i$ is equal to $\mathcal{C}_i \mathcal{C}_{n-2-i}$. The Catalan model is a branching Markov process where the splitting distribution is given by

$$\hat{q}_{n-1}(i) = \frac{\mathcal{C}_i \mathcal{C}_{n-2-i}}{\mathcal{C}_{n-1}}, \qquad i=0,\ldots,n-2.$$

The transformation $\Phi_n$ maps the Catalan model into the uniform model for PTs with $n$ leaves. The splitting distribution for PTs is then

$$q_n(i) = \frac{1}{2}\binom{n}{i}\frac{(2i-3)!!(2(n-i)-3)!!}{(2n-3)!!}, \qquad i = 1, \ldots, n-1.$$

This formula can also be found in [4].

LEMMA 2. *Let $q_n$ be a splitting distribution on $\{1, \ldots, n-1\}$. Let $h$ be a function of pairs of integers. Denote by $X_n$ an additive random variable, defined recursively as*

$$X_n \stackrel{d}{=} X_{I_n} + X_{J_n} + h(I_n, J_n),$$

*where the $I_\ell$'s are sampled under a branching Markov process of splitting distribution $q_n$ and where $I_n + J_n = n$. Define $\hat{X}_n$ by*

$$\hat{X}_n \stackrel{d}{=} \hat{X}_{\hat{I}_n} + \hat{X}_{\hat{J}_n} + h(\hat{I}_n + 1, \hat{J}_n + 1),$$

*where the $\hat{I}_\ell$'s are sampled under the branching Markov process of splitting distribution $\hat{q}_n$ with*

$$\hat{q}_n(i) = q_{n+1}(i+1), \qquad i = 0, \ldots, n-1,$$

*and where $\hat{I}_n + \hat{J}_n = n - 1$. Then we have*

$$X_n \stackrel{d}{=} \hat{X}_{n-1}.$$

PROOF. Note that $\hat{I}_n \stackrel{d}{=} I_{n+1} - 1$, that is, the distribution of $\hat{I}_n$ is given by $\hat{q}_n$. Similarly, we have

$$X_n \stackrel{d}{=} X_{\hat{I}_{n-1}+1} + X_{\hat{J}_{n-1}+1} + h(\hat{I}_{n-1}+1, \hat{J}_{n-1}+1).$$

Setting $\hat{X}_{n-1} = X_n$, we prove the result. □

REMARK 6. This lemma states that for additive random variables built from a Markov branching PT of splitting distribution $q_n$, there exist additive random variables $\hat{X}_n$ built from a Markov branching BST of distribution $\hat{q}_n$. In addition, $X_n$ and $\hat{X}_{n-1}$ have the same distribution. This lemma can obviously be generalized to multivariate random variables. In the next sections, all random variables are studied in the context of BSTs. Applying Lemma 2, the results can be transposed to PTs without difficulties.

**4. Yule model.** The Yule model is a branching Markov process for PTs with splitting probability

$$q_n(i) = \frac{1}{n-1} \quad \text{for } i = 1, \ldots, n-1.$$

Sackin's index $S_n$ has been defined as a sum of depths over the leaves. Sackin's index $S_n$ can also be expressed as a sum over the internal nodes [31],

$$S_n = \sum_{j=1}^{n-1} N_j,$$

where $N_j$ is the number of leaves descending from the internal node $j$. Applying Lemma 2, we obtain that $S_n$ has the same distribution as $\hat{S}_{n-1} + 2(n-1)$, where $\hat{S}_n$ is defined by

$$(8) \qquad \hat{S}_n \stackrel{d}{=} \hat{S}_{I_n} + \hat{S}'_{J_n} + n - 1,$$

$I_n$ is distributed uniformly over $\{0, \ldots, n-1\}$ and $J_n = n - 1 - I_n$. The recursion satisfied by $\hat{S}_n$ is well studied since it arises from the analysis of the quicksort algorithm or the internal path length of a BST under the random permutation model. Similarly, $C_n$ has the same distribution as $\hat{C}_{n-1}$, where

$$(9) \qquad \hat{C}_n \stackrel{d}{=} \hat{C}_{I_n} + \hat{C}'_{J_n} + |I_n - J_n|.$$

In order to describe the joint distribution of $(\hat{S}_n, \hat{C}_n)$ under the random permutation model, we shall follow the same lines of proof as Neininger [27] who studied the joint convergence of the Wiener index and the internal path length of a BST.

PROOF OF THEOREM 1. *Step* 1. *Computing expectations*. Denote $\hat{c}_n = \mathbb{E}[\hat{C}_n]$ and $\hat{s}_n = \mathbb{E}[\hat{S}_n]$. We have

$$(10) \qquad \hat{s}_n = 2n \log n + (2\gamma - 4)n + o(n).$$

We rewrite equation (9) as

$$(11) \qquad \hat{C}_n \stackrel{d}{=} \hat{C}_{I_n} + \hat{C}'_{J_n} + n - 1 - 2\min(I_n, J_n).$$

Conditioning on $I_n$ in the above equation, we find that

$$\hat{c}_n = (n - 1 - 2t_n) + \frac{2}{n} \sum_{k=0}^{n-1} \hat{c}_k,$$

where

$$t_n = \mathbb{E}[\min(I_n, J_n)] = \begin{cases} \dfrac{n-2}{4}, & \text{if } n \text{ is even,} \\ \dfrac{(n-1)^2}{4n}, & \text{if } n \text{ is odd.} \end{cases}$$

Applying Lemma 1 of [17], page 1691, we obtain that

$$\hat{c}_n = (n - 1 - 2t_n) + 2(n+1) \sum_{k=1}^{n-1} \frac{k - 1 - 2t_k}{(k+1)(k+2)}.$$

An asymptotic expansion of the above expression leads to the following result:

(12) $$\hat{c}_n = n \log n + (\gamma - 1 - \log 2)n + o(n).$$

*Step* 2. *Limit distribution.* Let us consider the rescaled quantities

$$\hat{X}_n = \begin{pmatrix} \dfrac{\hat{S}_n - \hat{s}_n}{n} \\ \dfrac{\hat{C}_n - \hat{c}_n}{n} \end{pmatrix}$$

and $\hat{X}'_n$, an independent copy of $\hat{X}_n$. From equations (8) and (11), we have

(13) $$\hat{X}_n = A_1^{(n)} \hat{X}_{I_n} + A_2^{(n)} \hat{X}'_{J_n} + b^{(n)},$$

where

$$A_1^{(n)} = \frac{1}{n} \begin{pmatrix} I_n & 0 \\ 0 & I_n \end{pmatrix},$$

$$A_2^{(n)} = \frac{1}{n} \begin{pmatrix} J_n & 0 \\ 0 & J_n \end{pmatrix}$$

and

$$\begin{pmatrix} b_S^{(n)} \\ b_C^{(n)} \end{pmatrix} = \frac{1}{n} \begin{pmatrix} \hat{s}_{I_n} + \hat{s}_{J_n} - \hat{s}_n + n - 1 \\ \hat{c}_{I_n} + \hat{c}_{J_n} - \hat{c}_n + n - 1 - 2\min(I_n, J_n) \end{pmatrix}.$$

Since $I_n/n$ converges in $L^2$ toward $U$, a uniform variable over $(0, 1)$, we have

$$A_1^{(n)} \xrightarrow{L^2} A_1^* = \begin{pmatrix} U & 0 \\ 0 & U \end{pmatrix}$$

and

$$A_2^{(n)} \xrightarrow{L^2} A_2^* = \begin{pmatrix} 1 - U & 0 \\ 0 & 1 - U \end{pmatrix}.$$

Using the asymptotic expansion of $\hat{c}_n$ given by (12) and the asymptotic expansion of $\hat{s}_n$ given by (10), we find that

$$b_S^{(n)} = \frac{1}{n}\left(2I_n \log\left(\frac{I_n}{n}\right) + 2J_n \log\left(\frac{J_n}{n}\right) + n\right) + o(1)$$

and

$$b_C^{(n)} = \frac{1}{n}\left(I_n \log\left(\frac{I_n}{n}\right) + J_n \log\left(\frac{J_n}{n}\right) + n - 2\min(I_n, J_n)\right) + o(1).$$

Thus, we have

$$\begin{pmatrix} b_S^{(n)} \\ b_C^{(n)} \end{pmatrix} \xrightarrow{L^2} b^* = \begin{pmatrix} 2U \log U + 2(1-U)\log(1-U) + 1 \\ U \log U + (1-U)\log(1-U) + 1 - 2\min(U, 1-U) \end{pmatrix}.$$

Assuming that $\hat{X}_n$ converges in distribution, the limiting distribution $\mathcal{L}(\hat{X})$ must satisfy the condition

(14) $$\hat{X} \stackrel{d}{=} A_1^* \hat{X} + A_2^* \hat{X}' + b^*,$$

where $\hat{X}$, $\hat{X}'$ and $(A_1^*, A_2^*, b^*)$ are independent and $\hat{X}' \stackrel{d}{=} \hat{X}$.

The multivariate contraction theorem [26] states that there is a unique probability distribution $\mathcal{L}(\hat{X})$ satisfying (14) in $\mathcal{M}_2$. Moreover, it states that the distribution of $\hat{X}_n$ converges toward the distribution of $\hat{X}$ in the Wasserstein–Mallow $d_2$-metric. The convergence in this metric is the same as the convergence in distribution and the convergence of second moments. Neininger's theorem can be applied provided that the following four conditions hold:

(i) $(A_1^{(n)}, A_2^{(n)}, b^{(n)}) \xrightarrow{L^2} (A_1^*, A_2^*, b^*), n \to \infty$,
(ii) $\mathbb{E}[\|(A_1^*)^t A_1^*\|_{\text{op}}] + \mathbb{E}[\|(A_2^*)^t A_2^*\|_{\text{op}}] < 1$,
(iii) $\mathbb{E}[\mathbb{1}_{\{I_n \le \ell\} \cup \{I_n = n\}}\|(A_1^*)^t A_1^*\|_{\text{op}}] \to 0$, for all $\ell \in \mathbb{N}, n \to \infty$,
(iv) $\mathbb{E}[\mathbb{1}_{\{J_n \le \ell\} \cup \{J_n = n\}}\|(A_2^*)^t A_2^*\|_{\text{op}}] \to 0$, for all $\ell \in \mathbb{N}, n \to \infty$,

where $\|A\|_{\text{op}} = \sup_{\|x\|=1} \|Ax\|$ is the operator norm of $A$. For symmetric matrices (which we are considering here) this equals the spectral radius.

(i) Has already been proved. (ii) Holds because

$$\mathbb{E}[\|(A_1^*)^t A_1^*\|_{\text{op}}] + \mathbb{E}[\|(A_2^*)^t A_2^*\|_{\text{op}}] = \mathbb{E}[U^2 + (1-U)^2] = \tfrac{2}{3} < 1.$$

(iii) and (iv) are obvious because $\|(A_r^*)^t A_r^*\|_{\text{op}} \le 1$ for $r = 1, 2$ and

$$\mathbb{P}(\{I_n \le \ell\} \cup \{I_n = n\}) = \mathbb{P}(\{J_n \le \ell\} \cup \{J_n = n\}) \le \frac{\ell+1}{n} \to 0$$

for all $\ell \in \mathbb{N}$ and $n \to \infty$. □

PROOF OF THEOREM 2. According to Theorem 1, equation (14) has a unique solution, so we can consider $(S, C)$ and $(S', C')$, two independent copies with $\mathcal{L}(S,C) = \mathcal{L}(S',C')$ being the fixed point of (14) and $U$ being uniform over $(0,1)$. By definition, we have

$$\mathcal{L}\begin{pmatrix} S \\ C \end{pmatrix} \stackrel{d}{=} \mathcal{TL}\begin{pmatrix} S \\ C \end{pmatrix}.$$

Using the fact that all random variables (except $U$) are centered, we find that

$$\mathbb{E}[C^2] = \mathbb{E}[U^2 C^2] + \mathbb{E}[(1-U)^2 C'^2]$$
$$+ \mathbb{E}[(U \log U + (1-U)\log(1-U) + 1 - 2\min(U, 1-U))^2].$$

Thus, we have

$$\text{Var}[C^2] = \mathbb{E}[C^2] = \left(3 - \frac{\pi^2}{6} - \log 2\right).$$

In the same way, we find that

$$\mathbb{E}[SC] = \mathbb{E}[U^2 SC] + \mathbb{E}[(1-U)^2 S'C']$$
$$+ \mathbb{E}[(2U \log U + 2(1-U)\log(1-U) + 1)$$
$$\times (U \log U + (1-U)\log(1-U) + 1 - 2\min(U, 1-U))].$$

This leads to

$$\text{Cov}(S, C) = \mathbb{E}[SC] = \frac{9}{2} - \frac{\pi^2}{3} - \log 2.$$

Using the fact that $\mathbb{E}[S^2] = 7 - 2\pi^2/3$ [32], we find that

$$\text{Cor}(S, C) = \frac{27 - 2\pi^2 - 6\log 2}{\sqrt{2(18 - \pi^2 - 6\log 2)(21 - 2\pi^2)}}.$$

Theorem 1 holds in the Wasserstein–Mallows $d_2$-metric, which implies the convergence of second moments. This leads to

$$\text{Var}[S_n] \sim \text{Var}[S] n^2, \qquad \text{Var}[C_n] \sim \text{Var}[C] n^2$$

and

$$\text{Cov}(S_n, C_n) \sim \text{Cov}(S, C) n^2, \qquad \text{Cor}(S_n, C_n) \sim \text{Cor}(S, C). \qquad \square$$

REMARK 7. Lemma 2 suggests that a more general class of statistics could be studied using the same technique. When the toll function $t_n = h(I_n, J_n)$ varies, general limit laws for recursive random variables can be found in [17]. Note that we have started the proof with the guess that the variance was of order $n$ (which may not be obvious in general). Readers interested in multivariate distributional recursion and convergence to a functional fixed-point solution could refer to the recent survey by Rüschendorf and Neininger [33].

**5. Uniform model.** For a given $n$, the uniform model assumes that all PTs with $n$ leaves are equally likely. Again, we use the fact that Sackin's and Colless' indices for a PT with $n$ leaves drawn according to the uniform model have the same probability distribution as $\hat{S}_{n-1} + 2(n-1)$ and $\hat{C}_{n-1}$, which are defined by equations (8) and (9), respectively. Under the Catalan model for BSTs, $I_n$ is distributed according to $\hat{q}_n$, where

$$\hat{q}_n(i) = \frac{\mathcal{C}_i \mathcal{C}_{n-1-i}}{\mathcal{C}_n}, \qquad i = 0, \ldots, n-1.$$

Conditional on $I_n$, $(\hat{S}_{I_n}, \hat{C}_{I_n})$ is independent of $(\hat{S}'_{J_n}, \hat{C}'_{J_n})$.

Clearly, $\hat{S}_n$ has the same distribution as the internal path length of a BST under the Catalan model and $\hat{C}_n$ has the same distribution as the random variable $\sum_{j=1}^{n-1} |\hat{L}_j - \hat{R}_j|$, where the sum is over the $n-1$ vertices of a BST drawn under the Catalan model. Note that $\hat{C}_n$ can be rewritten as

$$\sum_{j=1}^{n-1}(\hat{N}_j - 1) - 2\min(\hat{L}_j, \hat{R}_j),$$

where $\hat{N}_j$ is the number of vertices of the subtree rooted at $j$ (including $j$) and $\hat{L}_j$ ($\hat{R}_j$) is the number of vertices of the left (right) subtree. Then we have

$$\hat{C}_n = \hat{S}_n - 2\sum_{j=1}^{n-1} \min(\hat{L}_j, \hat{R}_j).$$

It is well known [37] that $\hat{S}_n/n^{3/2}$ converges in distribution to the Airy distribution. The proof relies on the one-to-one correspondence between binary trees and Bernoulli excursions. Takacs [37] computed the moments of $\hat{S}_n/n^{3/2}$ and their limiting values to establish convergence using the method of moments. The goal of this section is to prove the convergence

$$\mathbb{E}\left[\frac{\sum_{j=1}^{n-1} \min(\hat{L}_j, \hat{R}_j)}{n^{3/2}}\right] \to 0, \qquad n \to 0.$$

This implies that $(\hat{C}_n - \hat{S}_n)/n^{3/2}$ converges in probability to 0 and completes the proof of Theorem 3.

LEMMA 3. *Let $n \geq 2$ and consider a BST with $n$ vertices under the Catalan model. Denote by $j_0$ the root of the tree. We have*

$$\mathbb{E}[\min(\hat{L}_{j_0}, \hat{R}_{j_0})] \leq K\sqrt{n},$$

*for some constant $K$.*

PROOF. Stirling's formula yields a well-known asymptotic expansion for the Catalan number $\mathcal{C}_n$:

$$\mathcal{C}_n = \frac{4^n}{\sqrt{\pi n^3}}\left(1 + O\left(\frac{1}{n}\right)\right). \tag{15}$$

The expectation of $\min(\hat{L}_{j_0}, \hat{R}_{j_0})$ is given by

$$\mathbb{E}[\min(\hat{L}_{j_0}, \hat{R}_{j_0})] = \sum_{k=1}^{\lfloor n/2 \rfloor - 1} 2k \frac{\mathcal{C}_k \mathcal{C}_{n-1-k}}{\mathcal{C}_n} + \mathbb{1}_{\{n \in 2\mathbb{N}+1\}} \frac{(n-1)\mathcal{C}_{(n-1)/2}^2}{2\mathcal{C}_n}.$$

Using (15), we find that

$$\sum_{k=1}^{\lfloor n/2 \rfloor - 1} 2k \frac{\mathcal{C}_k \mathcal{C}_{n-1-k}}{\mathcal{C}_n} = \frac{4^{n-1}}{\pi \mathcal{C}_n (n-1)} \frac{1}{n-1}$$

$$\times \sum_{k=1}^{\lfloor n/2 \rfloor - 1} 2 \frac{(1 + O(1/k) + O(1/n))}{\sqrt{k/(n-1)}(1 - k/(n-1))^{3/2}}.$$

The sum in the right-hand side of the above equation is a Riemann sum. Using the fact that

$$\int_0^{1/2} \frac{1}{\sqrt{x}(1-x)^{3/2}} = 2,$$

we have

$$\sum_{k=1}^{\lfloor n/2 \rfloor - 1} 2k \frac{\mathcal{C}_k \mathcal{C}_{n-1-k}}{\mathcal{C}_n} = \frac{4^n}{\pi \mathcal{C}_n (n-1)}(1 + o(1)).$$

Using (15) again leads to

$$\mathbb{E}[\min(\hat{L}_{j_0}, \hat{R}_{j_0})] \sim \sqrt{n/\pi}, \tag{16}$$

which concludes the proof of the result. □

By conditioning on the sizes of the two subtrees of the root and using induction on $n$, it follows that

$$\mathbb{E}\left[\frac{\sum_{j=1}^{n-1} \min(\hat{L}_j, \hat{R}_j)}{n^{3/2}}\right] \leq K \mathbb{E}\left[\frac{\sum_j \sqrt{\hat{N}_j}}{n^{3/2}}\right].$$

In the following, we prove that the right-hand side of the above inequality converges to 0 as $n$ goes to $\infty$.

A lemma which is interesting in its own right provides the key argument. In the following, we use the standard convention that $\binom{0}{0} = 1$.

LEMMA 4. *Let $n \geq 1$. Consider a BST with $n$ vertices sampled according to the Catalan model. Pick a vertex $V$ at random from the $n$ vertices. Denoting $\hat{K}_n = \hat{N}_V$, we have*

$$\mathbb{P}(\hat{K}_n = k) = \frac{\mathcal{C}_k \binom{2n-2k}{n-k}}{n\mathcal{C}_n} \qquad \text{if } k = 1, \ldots, n.$$

PROOF. The proof relies on combinatorial arguments. Let us denote by $\nu_k(T)$ the number of subtrees with $k$ vertices in the BST $T$ having $|T| = n$ vertices. For $k = 1, \ldots, n$, we have

$$\mathbb{P}(\hat{K}_n = k) = \frac{1}{n} \sum_{j=1}^{n} \mathbb{P}(\hat{N}_j = k)$$

$$= \frac{1}{n} \mathbb{E}\left[\sum_{j=1}^{n} \mathbb{1}_{\{\hat{N}_j = k\}}\right] = \frac{1}{n} \mathbb{E}[\nu_k(T)].$$

$\nu_k(T)$ satisfies the linear recursion

$$(17) \qquad \nu_k(T) = \delta_{|T|,k} + \sum_{S} \nu_k(S),$$

where $\delta$ denotes the Kronecker symbol and the sum is over the subtrees of the root of $T$. Let us denote by $\mathcal{B}$ the set of all BSTs. We now introduce the cumulative generating function defined by

$$F_k(z) = \sum_{T \in \mathcal{B}} \nu_k(T) z^{|T|}$$

and

$$G_k(z) = \sum_{T \in \mathcal{B}} \delta_{|T|,k} z^{|T|} = \mathcal{C}_k z^k.$$

From the linear recurrence equation (17), Theorem 5.7 in [35] establishes the following relationship between cumulative generating functions $F_k$ and $G_k$:

$$F_k(z) = \frac{G_k(z)}{\sqrt{1-4z}}.$$

Using the fact that

$$\frac{1}{\sqrt{1-4z}} = \sum_{i \geq 0} \binom{2i}{i} z^i,$$

we find that

$$F_k(z) = \sum_{i \geq k} \binom{2i-2k}{i-k} \mathcal{C}_k z^i,$$

since the expectation of $\nu_k(z)$ is given by
$$\mathbb{E}[\nu_k(T)] = [z^n]F_k(z)/\mathcal{C}_n,$$
this completes the proof of the lemma. □

COROLLARY 1. *Let $n \geq 2$. Consider a PT with $n$ leaves sampled from the uniform model. Let $K_n$ denotes the number of leaves descending from a uniformly chosen random ancestor. We then have*
$$\mathbb{P}(K_n = k) = \frac{\mathcal{C}_{k-1}\binom{2n-2k}{n-k}}{(n-1)\mathcal{C}_{n-1}} \quad \text{for } k = 2, \ldots, n.$$

PROOF. This is a direct consequence of Lemma 2. □

REMARK 8. As $n$ goes to infinity, the distribution of $\hat{K}_n$ converges to
$$\mathbb{P}(\hat{K} = k) = 4^{-k}\mathcal{C}_k \sim \frac{1}{\sqrt{\pi}k^{3/2}}, \text{ for large } k.$$

The tail of the distribution of $\hat{K}$ has a power law with parameter $3/2$. This can be compared to a similar result in the context of BSTs [23] under the random permutation model. In this case, $\hat{K}_n$ has power law distribution with parameter 2. Since $3/2$ is less than 2, large random subtrees are more likely in the Catalan model than in the random permutation model. It was an expected result since Catalan binary trees are known to be more unbalanced than BSTs under the random permutation model (Section 5.6, [35]).

REMARK 9. Actually, the limiting distribution in the preceding comment is equal to the size of the critical Galton–Watson process with a binomial Bi(2, 1/2) offspring distribution (Lemma 9, [3]). This is the Galton–Watson process corresponding to binary trees.

We are now ready to prove that $\mathbb{E}[\sum_j \sqrt{\hat{N}_j}]/n^{3/2}$ converges to 0 as $n$ goes to $\infty$. Let $\alpha \in {]0,1[}$ and split the sum $\mathbb{E}[\sum_j \sqrt{\hat{N}_j}]$ into two parts:
$$\mathbb{E}\left[\sum_j \sqrt{\hat{N}_j}\right] = \mathbb{E}\left[\sum_{j,\hat{N}_j < n^\alpha} \sqrt{\hat{N}_j} + \sum_{j,\hat{N}_j \geq n^\alpha} \sqrt{\hat{N}_j}\right].$$

Obviously, we have
$$\frac{1}{n^{3/2}}\mathbb{E}\left[\sum_{j,\hat{N}_j < n^\alpha} \sqrt{\hat{N}_j}\right] \leq \frac{1}{n^{3/2}}\mathbb{E}\left[\sum_{j,\hat{N}_j < n^\alpha} \sqrt{n^\alpha}\right]$$
$$\leq n^{\alpha/2 - 1/2} \to 0$$

when $n$ goes to $\infty$. For the second term, we have

$$\frac{1}{n^{3/2}}\mathbb{E}\left[\sum_{j,\hat{N}_j\geq n^\alpha}\sqrt{\hat{N}_j}\right] \leq \frac{1}{n}\mathbb{E}\left[\sum_{j,\hat{N}_j\geq n^\alpha}1\right].$$

The right-hand side of the inequality is equal to $P(\hat{K}_n \geq n^\alpha)$. Applying Lemma 4, we find that

$$\frac{1}{n^{3/2}}\mathbb{E}\left[\sum_{j,\hat{N}_j\geq n^\alpha}1\right] \sim \kappa n^{-\alpha/2},$$

for some constant $\kappa$. This expression converges to 0 when $n$ goes to $\infty$. This completes the proof of Theorem 3.

REMARK 10. Fill and Kapur [12] established more precise results concerning $\mathbb{E}[\sum_j \sqrt{\hat{N}_j}]$. They proved that

(18) $$\mathbb{E}\left[\sum_j \sqrt{\hat{N}_j}\right] \sim \frac{1}{\sqrt{\pi}} n \log n.$$

Their results rely on Hadamard products. In a recent preprint, Ford [14] gave an alternate proof that $(\hat{C}_n - \hat{S}_n)/n^{3/2}$ converges in probability to 0. Note that our proof uses elementary arguments and is instructive in its own right as regards the shapes of Catalan trees. Besides, equation (18) follows easily from Lemma 4, by direct estimation of the sum by an integral.

PROOF OF THEOREM 4. In the following, we prove the convergence of the mixed moments of $\hat{S}_n$ and $\hat{C}_n$. For $k, l \geq 0$, we have

(19) $$\mathbb{E}[\hat{S}_n^k \hat{C}_n^\ell] \sim n^{3(k+\ell)/2}\mathbb{E}[\mathcal{A}^{k+\ell}].$$

The argument is similar to the argument given by Janson (Remark 3.5, [18]) to establish the convergence of the mixed moments of the internal path length and the Wiener index. Since the convergence in distribution has been established in Theorem 3, the above equation is equivalent to uniform integrability of $n^{-3(k+\ell)/2}\hat{S}_n^k \hat{C}_n^\ell$ for $n \geq 1$ and any fixed $k, \ell$. Since $\hat{C}_n \leq \hat{S}_n$, the result follows from the fact that $n^{-3k/2}(\hat{S}_n)^k$ is uniformly integrable for every fixed $k$. This is true because Takacs [37] proved the convergence of the moments of $\hat{S}_n$. $\square$

M. G. B. BLUM
O. FRANÇOIS
TIMC–TIMB
TEAM OF MATHEMATICAL BIOLOGY
PAVILLON TAILLEFER
FACULTÉ DE MÉDECINE
F38706 LA TRONCHE CEDEX
FRANCE
E-MAIL: michael.blum@imag.fr
olivier.francois@imag.fr

S. JANSON
DEPARTMENT OF MATHEMATICS
UPPSALA UNIVERSITY
P.O. BOX 480
S-751 06 UPPSALA
SWEDEN
E-MAIL: svante@math.uu.se